\def\ann{\mathop{\rm ann}\nolimits}
\def\ara{\mathop{\rm ara}\nolimits}
\def\Ass{\mathop{\rm Ass}\nolimits}

\def\Hom{\mathop{\rm Hom}\nolimits}
\def\InjH{\mathop{\rm E}\nolimits}
\def\LCMo{\mathop{\rm H}\nolimits}
\def\Min{\mathop{\rm Min}\nolimits}
\def\Spec{\mathop{\rm Spec}\nolimits}
\def\Supp{\mathop{\rm Supp}\nolimits}
\def\char{\mathop{\rm char}\nolimits}
\def\Naturalsign{{\rm l\kern-.23em N}}
\font \normal=cmr10 scaled \magstep0 \font \mittel=cmr10 scaled \magstep1 \font \gross=cmr10 scaled \magstep5
\input amssym.def
\input amssym.tex
\parindent=0pt
\gross On the associated primes of Matlis duals of top local cohomology
modules
\normal
\bigskip
Abstract:
\par
{\tt
After motivating the question we prove various results about the set of
associated primes
\par
of Matlis duals of
local cohomology modules. In easy cases we can calculate this set.
\par
An easy
application of this theory is the well-known fact that Krull dimension can be
ex-
\par
pressed by the vanishing of local cohomology modules.
}
\bigskip
\mittel
{\bf 0. Motivation and Notation}
\normal
\smallskip
Whenever $(R,\goth m)$ is a noetherian local ring, we denote by $\InjH
_R(R/\goth m)$ an $R$-injective hull of the residue field $R/\goth m$ and by
$D(M):=\Hom _R(M,\InjH _R(R/\goth m))$ the Matlis dual of $M$ (here $M$ is any
$R$-module, for basics on injective modules and Matlis duals see [3], [4],
[5], [9]).
\par
Let $I$ be an ideal of a noetherian local ring $(R,\goth m)$ such that $0=\LCMo ^2_I(R)=\LCMo
^3_I(R)=\dots $. For any $f\in I$ we have
$$\sqrt {I}=\sqrt {fR}\iff f \hbox { acts surjevtively on } \LCMo ^1_I(R)$$
(Proof: ``$\Rightarrow $'': Clear, because $f$ acts surjectively on $\LCMo
^1_{fR}(R)$. ``$\Leftarrow $'': For $f\in \goth p\in \Spec (R)$ the element
$f\in R$ acts like zero on $\LCMo ^1_I(R/\goth p)$ and at the same time
surjectively on $\LCMo ^1_I(R)\otimes _R(R/\goth p)=\LCMo ^1_I(R/\goth p)$; so
we must have $\LCMo ^1_I(R/\goth p)=0$ which implies $I\subseteq \goth p$.)
\smallskip
$f$ acts surjectively on $\LCMo ^1_I(R)$ if and only if $f$ acts
injectively on $D(\LCMo ^1_I(R))$, that is, if and only if $f$ is not
contained in any prime ideal associated to $D(\LCMo ^1_I(R))$.
\smallskip
The above statement is easily generalized to the following one: Let $I$ be an
ideal in a noetherian local ring $(R,\goth m)$, $n\in \Naturalsign $, $0=\LCMo
^{n+1}_I(R)=\LCMo ^{n+2}_I(R)=\dots $ and $f_1,\dots f_n\in I$. We assume
$I\neq R$ is a proper ideal to avoid trivial cases. Then
$$\eqalign {\sqrt {I}=\sqrt{(f_1,\dots ,f_n)R}&\iff f_i \hbox { acts
surjectively on }\LCMo ^{n+1-i}_I(R/(f_1,\dots ,f_{i-1})R)\hbox { for }
i=1,\dots ,n\cr &\iff f_i \hbox { acts
injectively on }D(\LCMo ^{n+1-i}_I(R/(f_1,\dots ,f_{i-1})R))\hbox { for
}i=1,\dots ,n\cr }$$
Under the additional assumption
$$\LCMo ^l_I(R)=0 \ \ \ \ (l\neq n)$$
(i. e. $f_1,\dots ,f_n$ form a regular sequence on $R$) we may formulate:
$$\sqrt {I}=\sqrt{(f_1,\dots ,f_n)R}\iff f_1,\dots ,f_n \hbox { is a regular
sequence on }D(\LCMo ^n_I(R))$$
This follows from
$$D(\LCMo ^n_I(R))/(f_1,\dots ,f_n)D(\LCMo ^n_I(R))=D(\Hom _R(R/(f_1,\dots
,f_n)R, \LCMo ^n_I(R)))\neq 0$$
and the fact that the short exact sequence
$$0\to R\buildrel f_n\over \to R\to R/f_nR\to 0$$
induces an exact sequence
$$0\to \LCMo ^{n-1}_I(R/f_nR)\to \LCMo ^n_I(R)\buildrel f_n\over \to \LCMo
^n_I(R)\to 0$$
and hence an equality
$$D(\LCMo ^{n-1}_I(R/f_nR))=D(\LCMo ^n_I(R))/f_nD(\LCMo ^n_I(R))$$
Thus the minimal number of generators of $I$ up to radical (this number is
often called the arithmetic rank $\ara (I)$ of $I$) is related to regular
sequences, an important concept in commutative algebra (basic facts on regular
sequences can be found in [5], [10], [11]).
\bigskip
\mittel
{\bf 1. Conjecture}
\normal
\smallskip
One may conjecture the following equality, being refered to by (*) from now
on:
\par
(*) Whenever $(R,\goth m)$ is a noetherian local ring,
$$\Ass _R(D(\LCMo ^i_{(x_1,\dots ,x_i)R}(R)))= \{ \goth p\in \Spec (R)\vert \LCMo ^i_{(x_1,\dots
,x_i)R}(R/\goth p)\neq 0\} $$
holds for any $i\geq 1$, $x_1,\dots ,x_i\in R$.
\bigskip
It is easy to see (paragraph 2.2) that we always have the inequality
$\subseteq $ in the above situation.
\bigskip
{\bf 1.1 Theorem}
\par
The following statements are equivalent:
\par
(i) Conjecture (*) holds.
\par
(ii) If $(R,\goth m)$ is an arbitrary noetherian local ring, $i\geq 1$ and
$x_1,\dots ,x_i\in R$, the set
$$Y:=\Ass _R(D(\LCMo ^i_{(x_1,\dots ,x_i)R}(R)))$$
has the following ``openness'' property (referred to as (O) ):
$\goth p,\goth q\in \Spec (R), \goth p\subseteq \goth q, \goth q\in
Y\Rightarrow \goth p\in Y$.
\par
(iii) Whenever $(R,\goth m)$ is a noetherian local domain, $i\geq 1$ and
$x_1,\dots ,x_i\in R$, $\LCMo ^i_{(x_1,\dots ,x_i)R}(R)\neq 0$ implies $\{ 0\}
\in \Ass _R(D(\LCMo ^i_{(x_1,\dots ,x_i)R}(R)))$
\par
Proof:
\par
(i) $\Rightarrow $(ii): In the situation of (ii) we have
$\Hom _R(R/\goth q,D(\LCMo ^i_{(x_1,\dots ,x_i)R}(R)))\neq 0$, hence
$$0\neq \Hom _R(R/\goth p,D(\LCMo ^i_{(x_1,\dots ,x_i)R})(R))=\Hom _R(\LCMo
^i_{(x_1,\dots ,x_i)R}(R/\goth p),\InjH _R(R/\goth m))$$
and so (i) implies $\goth p\in \Ass _R(D(\LCMo ^i_{(x_1,\dots
,x_i)R}(R)))$.
\par
(ii) $\Rightarrow $ (iii): $\LCMo ^i_{(x_1,\dots ,x_i)R}(R)\neq
0\Rightarrow D(\LCMo ^i_{(x_1,\dots x_i)R}(R))\neq 0\Rightarrow \Ass _R(D(\LCMo ^i_{(x_1,\dots
,x_i)R}(R)))\neq \emptyset \buildrel \hbox{(ii)}\over \Rightarrow \{ 0\} \in \Ass _R(D(\LCMo ^i_{(x_1,\dots
,x_i)R}(R)))$.
\par
(iii) $\Rightarrow $(i): Now $(R,\goth m)$ is an arbitrary
noetherian local ring and we have $\goth p\in \Spec (R)$ such that
$$\LCMo ^i_{(x_1,\dots x_i)R}(R/\goth p)\neq 0$$
We apply (iii) to the domain $R/\goth p$ (note $\LCMo ^i_{(x_1,\dots ,x_i)(R/\goth p)}(R/\goth p)=\LCMo
^i_{(x_1,\dots ,x_i)R}(R/\goth p)\neq 0$) and obtain an injection
$$\eqalign {R/\goth p&\to D(\LCMo ^i_{(x_1,\dots ,x_i)(R/\goth p)}(R/\goth
p))\cr &=\Hom _R(\LCMo ^i_{(x_1,\dots ,x_i)R}(R/\goth p),\InjH _R(R/\goth
m))\cr &=\Hom _R(\LCMo ^i_{(x_1,\dots ,x_i)R}(R)\otimes _R(R/\goth p), \InjH
_R(R/\goth m))\cr &=\Hom _R(R/\goth p,D(\LCMo ^i_{(x_1,\dots ,x_i)}(R)))\cr
&\subseteq D(\LCMo ^i_{(x_1,\dots ,x_i)R}(R))\cr }$$
proving (i).
\bigskip
{\bf 1.2 Remark}
\par
Suppose $(R,\goth m)$ is a noetherian local equicharacteristic domain,
$i\geq 1$ and $x_1,\dots ,x_i\in R$ are such that $\LCMo ^i_{(x_1,\dots
,x_i)R}(R)\neq 0$. Suppose furthermore one wants to show $\{ 0\} \in \Ass _R(D(\LCMo ^i_{(x_1,\dots
,x_i)R}(R)))$. For this one can replace $R$ by $\hat R/\goth p_0$, where $\hat
R$ is the $\goth m$-adic completion of $R$ and $\goth p_0\in \Spec (\hat R)$
is lying over the zero ideal of $R$, and so we may assume that $(R,\goth m)$
is a noetherian local equicharacteristic complete domain. Let $k$ be a
coefficient field of $R$. Now one can use a surjective $k$-algebra homomorphism
$k[[X_1,\dots ,X_n]]\to R$ (here $k[[X_1,\dots ,X_n]]$ denotes a power series
ring over $k$ in $n$ variables $X_1,\dots ,X_n$) mapping $X_1,\dots ,X_n$ to
$x_1,\dots x_n$, respectively, to reduce to the following problem:
\par
If $R=k[[X_1,\dots ,X_n]]$ is a power series ring over a field $k$ in $n$
variables $X_1,\dots ,X_n$, $1\leq i\leq n$, $\goth q\in \Ass _R(D(\LCMo
^i_{(X_1,\dots ,X_i)R}(R)))$, $\goth p\in \Spec (R)$, $\goth p\subseteq \goth
q$, then $\goth p\in \Ass _R(D(\LCMo ^i_{(X_1,\dots ,X_i)R}(R)))$ holds (that
is: The set $\Ass _R(D(\LCMo ^i_{(X_1,\dots ,X_i)R}(R)))$ has the above mentioned
``openness'' property (O) ).
\bigskip
\mittel
{\bf 2. Special Cases and their Consequences}
\normal
\smallskip
One can show that condition (iii) from theorem 1.1 (respectively the condition from remark 1.2) holds in special (regular)
situations; this is the content of the following paragraph:
\bigskip
\mittel
{\bf 2.1 Regular Cases}
\normal
\smallskip
{\bf 2.1.1 Lemma}
\par
Let $k$ be a field, $n\geq 1$, $R=k[[x_1,\dots ,x_n]]$ and $i\in \{
1,\dots ,n\} $. We set $I:=(x_1,\dots ,x_i)R$ and $\goth m:=(x_1,\dots
,x_n)R$. Then
$$\{ 0\} \in \Ass _R(D(\LCMo ^i_I(R)))$$
Proof:
\par
1. Case: $i=n$:
\par
Here $\LCMo ^i_I(R)=\InjH _R(R/\goth m)$ und also $D(\LCMo ^i_I(R))=R$ and the
statement follows.
\par
2. Case: $i<n$:
We have $\LCMo ^i_I(R)=\vtop{\baselineskip=1pt \lineskiplimit=0pt \lineskip=1pt\hbox{lim}
\hbox{$\longrightarrow $} \hbox{$^{^{l\in \bf N \setminus \{ 0\} }}$}}(R/(x_1^l,\dots
,x_i^l)R)$, the transition maps being induced by $R\to R$, $r\mapsto
(x_1\cdot \dots \cdot x_i)\cdot r$. So $D(\LCMo ^i_I(R))=\vtop{\baselineskip=1pt
\lineskiplimit=0pt \lineskip=1pt\hbox{lim}
\hbox{$\longleftarrow $} \hbox{$^{^{l\in \bf N\setminus \{ 0\} }}$}}(D(R/(x_1^l,\dots
,x_i^l)R))$; here $D(R/(x_1^l,\dots ,x_i^l)R)=\InjH _{R/(x_1^l,\dots
,x_i^l)R}(R/\goth m)$ ($\subseteq \InjH _R(R/\goth m)$), the transition maps
being induced by $\InjH _R(R/\goth m)\to \InjH _R(R/\goth m)$, $e\mapsto (x_1\cdot \dots \cdot
x_i)\cdot e$ and we have $\InjH _R(R/\goth m)=k[x_1^{-1},\dots
,x_n^{-1}]$. We define
$$\eqalign {\alpha :=&(1,x_1^{-1}\cdot \dots \cdot x_i^{-1}+x_{i+1}^{-1!}\cdot \dots
\cdot x_n^{-1!},\dots ,x_1^{-m}\cdot \dots \cdot x_i^{-m}+\cr &+(x_{i+1}^{-1!}\cdot
\dots \cdot x_n^{-1!})\cdot (x_1^{-(m-1)}\cdot \dots \cdot x_i^{-(m-1)})+\dots
+\cr &+(x_{i+1}^{-(m-1)!}\cdot \dots \cdot x_n^{-(m-1)!})\cdot (x_1^{-1}\cdot \dots
\cdot x_i^{-1})+x_{i+1}^{-m!}\cdot \dots \cdot x_n^{-m!},\dots )\in D(\LCMo
^i_I(R))\cr }$$
Here we consider the projective limit as a subset of a direct product. We
state $\ann _R(\alpha )=\{ 0\}$:
Assume there is an $f\in \ann _R(\alpha )\setminus \{ 0\} $. We choose
$(a_1,\dots ,a_n)\in \Supp (f)$ such that $(a_1,\dots ,a_i)$ is minimal
(using the ordering $(c_1,\dots c_i)\leq (c_1^\prime ,\dots ,c_i^\prime
):\iff c_1\leq c_1^\prime \wedge \dots \wedge c_i\leq c_i^\prime $) in $\{
(a_1^\prime ,\dots ,a_i^\prime )\vert \exists a_{i+1}^\prime ,\dots ,a_n^\prime
:(a_1^\prime ,\dots ,a_n^\prime )\in \Supp (f)\} $. We may assume $a_1=\max \{
a_1,\dots ,a_i\} $. We replace $f$ by $x_2^{a_1-a_2}\cdot \dots \cdot
x_i^{a_1-a_i}\cdot f$; this means $a_1=\dots =a_i=:a$. Choose $h_1,\dots ,h_i\in R$ and $g \in
k[[x_{i+1},\dots ,x_n]]\setminus \{ 0\} $ such that
$$f=x_1^{a+1}h_1+\dots +x_i^{a+1}h_i+(x_1^a\cdot \dots \cdot x_i^a)\cdot g$$
$f\cdot \alpha =0$ means: For every $m$ we have
$$\eqalign {0&=[ x_1^{a+1}h_1+\dots +x_i^{a+1}h_i+(x_1^a\cdot \dots \cdot x_i^a)\cdot
g] \cdot (x_1^{-m}\cdot \dots \cdot x_i^{-m}+\dots +x_{i+1}^{-m!}\cdot \dots
\cdot x_n^{-m!})\cr &=(x_1^{a+1}h_1+\dots +x_i^{a+1}h_i)\cdot [x_1^{-m}\cdot \dots
\cdot x_i^{-m}+\dots +(x_{i+1}^{-(m-a-1)!}\cdot \dots \cdot
x_n^{-(m-a-1)!})\cdot \cr &\cdot (x_1^{-(a+1)}\cdot \dots \cdot x_i^{-(a+1)})]+g\cdot
(x_1^{-(m-a)}\cdot \dots \cdot x_i^{-(m-a)}+\dots +x_{i+1}^{-(m-a)!}\cdot
\dots \cdot x_n^{-(m-a)!})\cr }$$
Choose $(b_{i+1},\dots ,b_n)$ minimal in $\Supp (g)$; then for all $m>>0$
the following statements must hold:
$$\eqalign {(m-a)!-b_{i+1}&\leq (m-a-1)!\cr &\vdots \cr (m-a)!-b_n&\leq
(m-a-1)!\cr }$$
For $m>>0$ this leads to a contradiction, the assumption is wrong and the
lemma is proved.
\bigskip
{\bf 2.1.2 Lemma}
\par
Let $p$ be a prime number, $C$ a complete $p$-ring, $n\geq 1$, $R=C[[x_1,\dots
,x_n]]$ and $i\in \{ 1,\dots ,n\} $. We set $I:=(x_1,\dots ,x_i)R$ and $\goth
m:=(p,x_1,\dots ,x_n)R$. Then $\{ 0\} \in \Ass _R(D(\LCMo ^i_I(R)))$.
\par
Proof:
\par
We have $\LCMo ^i_I(R)=\vtop{\baselineskip=1pt \lineskiplimit=0pt \lineskip=1pt\hbox{lim}
\hbox{$\longrightarrow $} \hbox{$^{^{l\in \bf N \setminus \{ 0\}
}}$}}(R/(x_1^l,\dots ,x_i^l)R)$, the transition maps being induced by $R\to R$,
$r\mapsto (x_1\cdot \dots \cdot x_i)\cdot r$. We deduce $D(\LCMo ^i_I(R))=\vtop{\baselineskip=1pt
\lineskiplimit=0pt \lineskip=1pt\hbox{lim}
\hbox{$\longleftarrow $} \hbox{$^{^{l\in \bf N\setminus \{ 0\}
}}$}}(D(R/(x_1^l,\dots ,x_i^l)R))$; we recall $D(R/(x_1,\dots ,x_i^l)R)=\InjH
_{R/(x_1^l,\dots x_i^l)R}(R/\goth m)$ ($\subseteq \InjH _R(R/\goth m)$), the
transition maps being induced by $\InjH _R(R/\goth m)\to \InjH _R(R/\goth m)$,
$e\mapsto (x_1\cdot \dots \cdot x_i)\cdot e$. Furthermore
$$\InjH _R(R/\goth m)=(C_p/C)[x_1^{-1},\dots ,x_n^{-1}]$$ holds (because of
$\InjH _R(R/\goth m)=\LCMo ^{n+1}_{(p,x_1,\dots ,x_n)R}(R)=\LCMo
^1_{pR}(R)\otimes _R\dots \otimes _R\LCMo ^1_{x_nR}(R)=(C_p/C)\otimes
_C((R_{x_1}/R)\otimes _R\dots \otimes _R(R_{x_n}/R))$).
We define
$$\eqalign {\alpha :=&(p^{-1},p^{-1}x_1^{-1}\cdot \dots \cdot
x_i^{-1}+p^{-1!}x_{i+1}^{-1!}\cdot \dots \cdot x_n^{-1!},\dots
,p^{-1}x_1^{-m}\cdot \dots \cdot x_i^{-m}+\cr &+(p^{-1!}x_{i+1}^{-1!}\cdot \dots
\cdot x_n^{-1!})\cdot (x_1^{-(m-1)}\cdot \dots \cdot x_i^{-(m-1)})+\dots
+\cr &+(p^{-(m-1)!}x_{i+1}^{-(m-1)!}\cdot \dots \cdot x_n^{-(m-1)!})\cdot
(x_1^{-1}\cdot \dots \cdot x_i^{-1})+p^{-m!}x_{i+1}^{-m!}\cdot \dots \cdot
x_n^{-m!},\dots )\in D(\LCMo ^i_I(R))\cr }$$
and, similar to the proof of lemma 2.1.1, we show that $\ann _R(\alpha
)=0$. Assume to the contrary there is an $f\in \ann _R(\alpha )\setminus \{
0\} $. Choose $(a_1,\dots ,a_i)$ minimal in
$$\{ (a_1^\prime, \dots ,a_i^\prime )\vert \hbox{there exists }a_{i+1}^\prime
,\dots a_n^\prime \hbox{ such that }(a_1^\prime ,a_n^\prime )\in \Supp (f)\}
$$
Like before we may assume $a_1=\dots =a_i=:a$. Choose $h_1,\dots ,h_i\in R$
and $g\in C[[x_{i+1},\dots ,x_n]]\setminus \{ 0\} $ such that
$$f=x_1^{a+1}\cdot h_1+\dots +x_i^{a+1}\cdot h_i+x_1^a\cdot \dots \cdot
x_i^a\cdot g$$
$\alpha \cdot f=0$ implies, for all $m\in \bf N \setminus \{ 0\} $,
$$\eqalign {0&=(x_1^{a+1}h_1+\dots +x_i^{a+1}h_i+x_1^a\cdot \dots \cdot x_i^a\cdot
g)\cdot (p^{-1}x_1^{-m}\cdot \dots \cdot x_i^{-m}+\dots
+p^{-m!}x_{i+1}^{-m!}\cdot \dots \cdot x_n^{-m!})\cr &=(x_1^{a+1}h_1+\dots
+x_i^{a+1}h_i)\cdot [p^{-1}x_1^{-m}\cdot \dots \cdot x_i^{-m}+\dots
+(p^{-(m-a-1)!}\cdot x_{i+1}^{-(m-a-1)!}\cdot \dots \cdot
x_n^{-(m-a-1)!})\cdot \cr &\cdot (x_1^{-(a+1)}\cdot \dots \cdot x_i^{-(a+1)})]+g\cdot
(p^{-1}x_1^{-(m-a)}\cdot \dots \cdot x_i^{-(m-a)}+\dots
+p^{-(m-a)!}x_{i+1}^{-(m-a)!}\cdot \dots \cdot x_n^{-(m-a)!})}$$
Let $(b_{i+1},\dots ,b_n)$ be minimal in $\Supp (g)$ and $c\in C$ be the
coefficient of $g$ in front of $x_{i+1}^{b_{i+1}}\cdot \dots \cdot
x_n^{b_n}$. In $C_p/C$ we have $c\cdot p^{-(m-a)!}\neq 0$ for all $m>>0$. So,
like before, we must have
$$(m-a)!-b_{i+1}\leq (m-a-1)!$$
$$(m-a)!-b_n\leq (m-a-1)!$$
for all $m>>0$, which leads to a contradiction again.
\bigskip
{\bf 2.1.3 Lemma}
\par
Let $p$ be a prime number, $C$ a complete $p$-ring, $n\in \bf N$,
$R=C[[x_1,\dots ,x_n]]$, $i\in \{ 0,\dots ,n\} $, $I:=(p,x_1,\dots ,x_i)R$ and
$\goth m:=(p,x_1,\dots ,x_n)R$. Then
$$\{ 0\} \in \Ass _R(D(\LCMo^{i+1}_I(R)))$$
Proof:
\par
1. Case: $i=n$: In this case we have $\LCMo ^{i+1}_I(R)=\InjH _R(R/\goth m)$
   and hence $D(\LCMo ^{i+1}_I(R))=R$.
\par
2. Case: $i<n$: Similar to the situation in the proof of lemma 2.1.2 we have
$$D(\LCMo ^{i+1}_I(R))=\vtop{\baselineskip=1pt
\lineskiplimit=0pt \lineskip=1pt\hbox{lim}
\hbox{$\longleftarrow $}}(\InjH _{R/(p,x_1,\dots ,x_i)R}(R/\goth m)\buildrel
p\cdot x_1\cdot \dots \cdot x_i\over \longleftarrow \InjH _{R/(p^2,x_1^2,\dots ,x_i^2)R}(R/\goth m)\buildrel
p\cdot x_1\cdot \dots \cdot x_i\over \longleftarrow \dots )$$
$$\InjH _R(R/\goth m)=(C_p/C)[x_1^{-1},\dots ,x_n^{-1}]$$
and we define
$$\eqalign {\alpha :=&(p^{-1},p^{-2}x_1^{-1}\cdot \dots \cdot
x_i^{-1}+p^{-2}x_{i+1}^{-1!}\cdot \dots \cdot x_n^{-1!},\dots
,p^{-(m+1)}x_1^{-m}\cdot \dots \cdot x_i^{-m}+\cr &+p^{-(m+1)}x_{i+1}^{-1!}\cdot
\dots \cdot x_n^{-1!}\cdot x_1^{-(m+1)}\cdot \dots \cdot x_i^{-(m-1)}+\dots
+\cr &+p^{-(m+1)}x_{i+1}^{-(m-1)!}\cdot \dots \cdot x_n^{-(m-1)!}\cdot x_1^{-1}\cdot
\dots \cdot x_i^{-1}+p^{-(m+1)}x_{i+1}^{-m!}\cdot \dots \cdot x_n^{-m!},\dots
)\cr }$$
Again we state $\ann _R(\alpha )$ and assume, to the contrary, that there
exists an $f\in \ann _R(\alpha )\setminus \{ 0\} $, choose $(a_1,\dots ,a_i)$
minimal in
$$\{ (a_1^\prime ,\dots ,a_i^\prime )\vert \hbox{There exist }a_{i+1}^\prime
,\dots ,a_n^\prime \hbox{ such that }(a_1^\prime ,\dots ,a_n^\prime )\in \Supp
(f)\} $$
may assume $a_1=\dots =a_i=:a$ and choose $h_1,\dots ,h_i\in R$, $g\in
C[[x_{i+1},\dots ,x_n]]$ such that
$$f=x_1^{a+1}h_1+\dots +x_i^{a+1}h_i+(x_1^a\cdot \dots \cdot a_i^a)\cdot g$$
This means, for all $m\in \bf N$,
$$\eqalign {0&=(x_1^{a+1}h_1+\dots +x_i^{a+1}h_i)[p^{-(m+1)}x_1^{-m}\cdot \dots \cdot
x_i^{-m}+\dots +(p^{-(m+1)}x_{i+1}^{-(m-a-1)!}\cdot \dots \cdot
x_n^{-(m-a-1)!})\cdot \cr &\cdot (x_1^{-(a+1)}\cdot \dots \cdot x_i^{-(a+1)})]+g\cdot
(p^{-(m+1)}x_1^{-(m-a)}\cdot \dots \cdot x_i^{-(m-a)}+\dots
+p^{-(m+1)}x_{i+1}^{-(m-a)!}\cdot \dots \cdot x_n^{-(m-a)!})\cr }$$
Choose $(b_{i+1},\dots ,b_n)$ minimal in $\Supp (g)$ and let $c\in C$ be the
coefficient of $g$ in front of $x_{i+1}^{b_{i+1}}\cdot \dots \cdot
x_n^{b_n}$. In $C_p/C$ we have $g\cdot p^{-(m+1)}\neq 0$ for all $m>>0$, and
so we must have for all $m>>0$
$$(m-a)!-b_{i+1}\leq (m-a-1)!$$
$$(m-a)!-b_n\leq (m-a-1)!$$
which leads to a contradiction, proving the lemma.
\bigskip
\mittel
{\bf 2.2 Some special Associated Primes of Duals of Top Local Cohomology Modules}
\normal
\smallskip
{\bf 2.2.1 Theorem}
\par
Let $(R,\goth m)$ be a noetherian local ring, $i\geq 1$ and $x_1,\dots ,x_i\in R$. Then
$$\Ass _R(D(\LCMo ^i_{(x_1,\dots ,x_i)R}(R)))\subseteq \{ \goth p\in \Spec (R)\vert \LCMo ^i_{(x_1,\dots
,x_i)R}(R/\goth p)\neq 0\} $$
holds. If $R$ is equicharacteristic,
$$\{ \goth p\in \Spec (R)\vert x_1,\dots ,x_i\hbox { is part of a
system of parameters for }R/\goth p\} \subseteq \Ass _R(D(\LCMo ^i_{(x_1,\dots
,x_i)R}(R)))$$
holds.
\par
Proof:
\par
For the first inclusion let $\goth p$ be an arbitrary element of
$\Ass _R(D(\LCMo ^i_{(x_1,\dots ,x_i)R}(R)))$; in particular
$$\eqalign {0&\neq \Hom _R(R/\goth p,D(\LCMo ^i_{(x_1,\dots ,x_i)R}(R)))\cr
&=\Hom _R(R/\goth p,\Hom _R(\LCMo ^i_{(x_1,\dots ,x_i)R}(R),\InjH _R(R/\goth
m)))\cr &=\Hom _{R/\goth p}(\LCMo ^i_{(x_1,\dots ,x_i)R}(R)\otimes _RR/\goth
p,\InjH _{R/\goth p}(R/\goth m))\cr &=\Hom _{R/\goth p}(\LCMo ^i_{(x_1,\dots
,x_i)R}(R/\goth p),\InjH _{R/\goth p}(R/\goth m))\cr }$$
implying $\LCMo ^i_{(x_1,\dots ,x_i)R}(R/\goth p)\neq 0$.
\par
For the second inclusion let $\goth p\in \Spec (R)$ and $x_{i+1},\dots ,x_n\in
R$ such that $x_1,\dots ,x_n$ (more precisely: their images in $R/\goth p$) is
a system of parameters for $R/\goth p$; then $n=\dim (R/\goth p)$. $x_1,\dots
,x_n$ is a system of parameters also in $\hat R/\goth p\hat R$. Choose $\goth
q\in \Spec (\hat R)$ with $\dim (\hat R/\goth q)=\dim (R/\goth p)$. This
implies $\goth q\in \Min (\hat R)$ and $\goth q\cap R=\goth p$. Because of
$\dim (\hat R/\goth q)=\dim (R/\goth p)$ the elements $x_1,\dots
,x_n$ form a system of parameters of $\hat R/\goth q$. It is sufficient to show $ \goth q\in \Ass
_{\hat R}(D(\LCMo ^i_{(x_1,\dots ,x_i)\hat R}(\hat R)))$, because we have
$$\eqalign {D(\LCMo ^i_{(x_1,\dots ,x_i)\hat R}(\hat R))&=\Hom _{\hat R}(\LCMo
^i_{(x_1,\dots ,x_i)\hat R}(\hat R),\InjH _{\hat R}(\hat R/\goth m\hat R))\cr &=\Hom _{\hat R}(\LCMo
^i_{(x_1,\dots ,x_i)\hat R}(\hat R),\InjH _R(R/\goth m))\cr &=\Hom _{\hat
R}(\LCMo ^i_{(x_1,\dots ,x_i)R}(R)\otimes _R\hat R,\InjH _R(R/\goth m))\cr
&=\Hom _R(\LCMo ^i_{(x_1,\dots ,x_i)R}(R),\InjH _R(R/\goth m))\cr &=D(\LCMo
^i_{(x_1,\dots ,x_i)R}(R))\cr }$$
and so every monomorphism $\hat R/\goth q\to D(\LCMo
^i_{(x_1,\dots ,x_i)\hat R}(\hat R))$ induces a monomorphism
$$R/\goth p\buildrel \hbox {kan.}\over \to \hat R/\goth q\to D(\LCMo ^i_{(x_1,\dots ,x_i)\hat
R}(\hat R))=D(\LCMo ^i_{(x_1,\dots ,x_i)R}(R))$$
Put into different words this means we may assume that $R$ is complete.
\par
We have to show that the zero ideal of $R/\goth p$ is associated to $$\Hom _R(R/\goth p,D(\LCMo^i_{x_1,\dots ,x_i)R}(R))=
D(\LCMo ^i_{(x_1,\dots ,x_i)R/\goth p}(R/\goth p))$$
(this equality was shown in the proof of the second inclusion). Replacing $R$ by $R/\goth p$
we may assume that $R$ is a domain and $\goth p$ is the zero ideal in $R$. Let
$k\subseteq R$ denote a coefficient field.
$$R_0:=k[[x_1,\dots ,x_n]]\subseteq R$$
is an $n$-dimensional regular local subring of $R$, over which $R$ is
module-finite. Let $\goth m_0$ denote the maximal ideal of $R_0$. The $R$-Modul
$\Hom _{R_0}(R,\InjH _{R_0}(R_0/\goth m_0))$ is isomorphic to $\InjH
_R(R/\goth m)$. We have
$$\eqalign {D(\LCMo ^i_{(x_1,\dots ,x_i)R}(R))&=
\Hom _R(\LCMo ^i_{(x_1,\dots ,x_i)R}(R),\InjH _R(R/\goth m))\cr
&=\Hom _R(\LCMo ^i_{(x_1,\dots ,x_i)R_0}(R_0)\otimes
_{R_0}R,\InjH _R(R/\goth m))\cr
&=\Hom _{R_0}(\LCMo ^i_{(x_1,\dots ,x_i)R_0}(R_0),
\Hom _{R_0}(R,\InjH _{R_0}(R_0/\goth m_0)))\cr
&=\Hom _{R_0}(R,\Hom _{R_0}(\LCMo ^i_{(x_1,\dots ,x_i)R_0}
(R_0),\InjH _{R_0}(R_0/\goth m_0)))\cr
&=\Hom _{R_0}(R,D(\LCMo ^i_{(x_1,\dots ,x_i)R_0}(R_0)))\cr }$$
by lemma 2.1.1 there exists a monomorphism $R_0\to
D(\LCMo ^i_{(x_1,\dots ,x_i)R_0}(R_0))$; it induces a monomorphism
$$\Hom _{R_0}(R,R_0)\to \Hom _{R_0}(R,D(\LCMo ^i
_{(x_1,\dots ,x_i)R_0}(R_0)))=D(\LCMo ^i_{(x_1,\dots ,x_i)R}(R))$$
$R$ is a domain and module-finite over $R_0$, and thus
$\{ 0\} \in \Supp _R(\Hom _{R_0}(R,R_0))$; the statement follows now.
\bigskip
In case $i=1$ theorem 2.2.1 means:
\bigskip
{\bf 2.2.2 Corollary}
\par
Let $(R,\goth m)$ be a noetherian local equicharacteristic ring and $x\in R$. Then
$$\Ass _R(D(\LCMo ^1_{xR}(R)))=\Spec (R)\setminus \goth V(x)$$
\bigskip
{\bf 2.2.3 Theorem}
\par
Let $(R,\goth m)$ be a noetherian local ring of mixed characteristic, $p=\char
(R/\goth m)$, $i\geq 0$ and $x_1,\dots ,x_i\in R$. Then
$$\{ \goth p\in \Spec (R)\vert p,x_1,\dots ,x_i\hbox { is part of a system of
parameters for }R/\goth p\} \subseteq \Ass _R(D(\LCMo ^{i+1}_{(p,x_1,\dots
,x_i)R}(R)))$$
In case $i\geq 1$
$$\{ \goth p\in \Spec (R)\vert p,x_1,\dots ,x_i\hbox { is part of a system of
parameters for }R/\goth p\} \subseteq \Ass _R(D(\LCMo ^i_{(x_1,\dots
,x_i)R}(R)))$$
holds in addition.
\par
Theorem 2.2.3 is proved in a similar way like Theorem 2.2.1, using lemma 2.1.2 and lemma 2.1.3
instead of lemma 2.1.1.
\bigskip
{\bf 2.2.4 Remarks}
\par
(i) If one has $\Ass _R(D(\LCMo ^i
_{(x_1,\dots ,x_i)R}(R)))=\emptyset $ in the situation of the theorem, it
follows that $\LCMo ^i_{(x_1,\dots ,x_i)R}(R)=0$ and also $\LCMo ^i
_{(x_1,\dots ,x_i)R}(R/\goth p)=0$ for every $\goth p
\in \Spec (R)$ (by a well-known theorem).
\smallskip
(ii) The second inclusion of theorem 2.2.1 is not an equality in general: For a
counterexample let $k$ be a field, $R=k[[y_1,y_2,y_3,y_4]]$ and
define $x_1=y_1y_3$, $x_2=y_2y_4$, $x_3=y_1y_4+y_2y_3$. $x_1,x_2,x_3$ is not part of a system of
parameters for $R$, but we have $\sqrt {(x_1,x_2,x_3)R}=(y_1,y_2)R\cap
(y_3,y_4)R$ and so a Mayer-Vietoris sequence argument shows $\InjH _R(k)=\LCMo
^3_{(y_1,y_2)R\cap (y_3,y_4)R}(R)=\LCMo ^3_{(x_1,x_2,x_3)R}(R)$ and so $D(\LCMo
^3_{(x_1,x_2,x_3)R}(R))=R$ implying $\{ 0\} \in \Ass _R(D(\LCMo
^3_{(x_1,x_2,x_3)R}(R)))$.
\smallskip
(iii) In the situation of theorem 2.2.1 let denote
$$Z_1:=\{ \goth p\in \Spec (R)\vert \LCMo ^i_{(x_1,\dots ,x_i)R}(R/\goth
p)\neq 0\} $$
$$Z_2:=\{ \goth p\in \Spec (R)\vert x_1,\dots ,x_i\hbox { is part of a
system of parameters for }R/\goth p\}$$
Then $Z_1$ has property (O) and, if we assume in addition that $R$ is regular,
$Z_2$ also has (O).
\bigskip
\mittel
{\bf 2.3 Comparison between two Matlis Duals}
\normal
\smallskip
Let $(R,\goth m)$ be a noetherian local equicharacteristic complete ring with
coefficient field $k$,
$y_1,\dots ,y_i\in R$ such that $R_0:=k[[y_1,\dots ,y_i]]$ is regular and of
dimension $i$ (e. g. when $\LCMo ^i_{(y_1,\dots ,y_i)R}(R)\neq 0$). Let $D_R$
denote the Matlis-dual functor with respect to $R$ and $D_{R_0}$ the one with
respect to $R_0$. By local duality we get
$$D_{R_0}(\LCMo ^i_{(y_1,\dots ,y_i)R}(R))=\Hom _{R_0}(R\otimes _{R_0}\LCMo
^i_{(y_1,\dots ,y_i)R_0}(R_0),\InjH _{R_0}(k))=\Hom _{R_0}(R,R_0)$$
$\Hom _{R_0}(R,\InjH _{R_0}(k))$ is an injective $R$-module; there exists an
injective $R$-module $E^\prime $ such that
$$\Hom _{R_0}(R,\InjH _{R_0}(k))=\InjH _R(k)\oplus E^\prime $$
We set $E=\Gamma _{(y_1,\dots ,y_i)R}(E^\prime )$. An easy computation shows
$$D_{R_0}(\LCMo ^i_{(y_1,\dots ,y_i)R}(R))=D_R(\LCMo ^i_{(y_1,\dots
,y_i)R}(R))\oplus \Hom _R(\LCMo ^i_{(y_1,\dots ,y_i)R}(R),E)$$
and hence
$$\Ass _R(D_{R_0}(\LCMo ^i_{(y_1,\dots ,y_i)R}(R)))=\Ass _R(D_R(\LCMo ^i_{(y_1,\dots
,y_i)R}(R)))\cup \Ass _R(\Hom _R(\LCMo ^i_{(y_1,\dots ,y_i)R}(R),E))$$
It is natural to ask for relations between $D_R(\LCMo ^i_{(y_1,\dots
,y_i)R}(R))$ and $D_{R_0}(\LCMo ^i_{(y_1,\dots
,y_i)R}(R))$.
\bigskip
For every $\goth p\in Z:=\{ \goth p\in \Spec (R)\vert (y_1,\dots ,y_i)R\subseteq
\goth p\subsetneq \goth m\} $ we choose a set $\mu_\goth p$ such that
$$E=\bigoplus _{\goth p\in Z}\InjH _R(R/\goth p)^{(\mu_\goth p)}$$
\bigskip
{\bf 2.3.1 Remark}
\par
One has $\mu_\goth p\neq \emptyset $ for every $ \goth p\in Z$.
\par
Proof:
\par
We have to show that $\goth p$ is associated to $\Hom _{R_0}(R/\goth p,\InjH
_{R_0}(k))$; the last module is $\Hom _{R_0}(R/\goth p,k)$, because $\goth p$
is annihilated by $y_1,\dots ,y_i$; that means we have to show that if
$(R,\goth m)$ is a noetherian local equicharacteristic complete domain with
coefficient field $k$, then the zero ideal is associated to $\Hom _k(R,k)$:
\par
Let $x_1,\dots ,x_n\in R$ be a system of parameters for $R$. Then
$R_0:=k[[x_1,\dots ,x_n]]$ is a regular subring of $R$, over which $R$ is
module-finite. One has $\Hom _k(R,k)=\Hom _{R_0}(R,\Hom _k(R_0,k))$ and
therefore it is sufficient to prove $\{ 0\} \in \Ass _{R_0}(\Hom _k(R_0,k))$,
because then every $R_0$-injection
$$R_0\to \Hom _k(R_0,k)$$
induces an $R$-injection
$$\Hom _{R_0}(R,R_0)\to \Hom _k(R,k)$$
and one always has $\{ 0\} \in \Supp _R(\Hom _{R_0}(R,R_0))$; thus we may
assume $R=k[[x_1,\dots ,x_n]]$ from now on:
\par
For $i=1,\dots ,n$ we set $R_i:=k[[x_1,\dots ,x_i]]$. Again we have
$$\Hom _k(R_i,k)=\Hom _{R_{i-1}}(R_i,\Hom _k(R_{i-1},k))$$
for $i=2,\dots ,n$. Using this and an obvious induction argument, the statement follows from lemma
2.3.2.
\bigskip
{\bf 2.3.2 Lemma}
\par
Let $k$ be a field and $R_0:=k[[x_1,\dots ,x_n]]$, $R:=k[[x_1,\dots
,x_n,x]]=R_0[[x]]$ be power series rings in the variables $x_1,\dots ,x_n,x$,
respectively. Then
$$\{ 0\} \in \Ass _R(\Hom _{R_0}(R,R_0))$$
holds.
\par
Proof:
\par
By $\goth m_0$ we denote the maximal ideal of $R_0$.
The canonical short exact sequence
$$0\to R_0[x]\to R_0[[x]]\to R_0[[x]]/R_0[x]\to 0$$
induces an exact sequence
$$0\to \Hom _{R_0}(R_0[[x]]/R_0[x],R_0)\to \Hom _{R_0}(R_0[[x]],R_0)\buildrel
\alpha \over \to \Hom _{R_0}(R_0[x],R_0)$$
The map $\alpha $ is the dual of the canonical map
$$\LCMo ^n_{\goth m_0}(R_0[x])=\LCMo ^n_{\goth
m_0}(R_0)\otimes _{R_0}R_0[x]\to \LCMo ^n_{\goth m_0}(R_0)\otimes
R_0[[x]]=\LCMo ^n_{\goth m_0}(R_0[[x]])$$
which is obviously injective. This means that $\alpha $ is surjective.
The $R_0[x]$-module $\Hom _{R_0}(R_0[x],R_0)$ can be written as
$R_0[[x^{-1}]]$ and in $R_0[[x^{-1}]]$ the element $h^\prime
:=1+x^{-1!}+x^{-2!}+\dots $ has $R_0[x]$-annihilator zero. Choose an element
$h\in \Hom _{R_0}(R_0[[x]],R_0)$ which is mapped to $h^\prime $ by $\alpha $;
then $\ann _{R_0[x]}(h)=\{ 0\} $ which implies $\ann _{R_0[[x]]}(h)=\{ 0\} $.
\bigskip
Quite generally every $\InjH _R(R/\goth p)$ is an $R_\goth p$-module and so we
always have
$$\Hom _R(\LCMo ^i_{(y_1,\dots ,y_i)R}(R),\InjH _R(R/\goth p))=\Hom _{R_\goth
p}(\LCMo ^i_{(y_1,\dots ,y_i)R_\goth p}(R_\goth p),E_{R_\goth p}(R_\goth
p/\goth pR_\goth p))$$
\bigskip
{\bf 2.3.3 Remark}
\par
If conjecture (*) holds, there are the following consequences:
\par
(i)
$$\Ass _R(\Hom _R(\LCMo ^i_{(y_1,\dots ,y_i)R}(R),\InjH _R(R/\goth p)))=\{
\goth q\in \Spec (R)\vert (\LCMo ^i_{(y_1,\dots ,y_i)R}(R/\goth q))_\goth
p\neq 0\} $$
In particular we have
$$\Ass _R(\Hom _R(\LCMo ^i_{(y_1,\dots ,y_i)R}(R),\InjH _R(R/\goth
p)))\subseteq \Ass _R(\Hom _R(\LCMo ^i_{(y_1,\dots ,y_i)R}(R),\InjH
_R(k)))$$
(ii) We have
$$\eqalign {\bigoplus _{\goth p\in Z}\Hom _R(\LCMo ^i_{(y_1,\dots ,y_i)R}(R),\InjH
_R(R/\goth p))^{(\mu_\goth p)}&\subseteq \Hom _R(\LCMo ^i_{(y_1,\dots
,y_i)R}(R),E)\cr &\subseteq \prod _{\goth p\in Z}\Hom _R(\LCMo ^i_{(y_1,\dots ,y_i)R}(R),\InjH
_R(R/\goth p))^{\mu_\goth p}\cr }$$
Every non-trivial annulator of an element of $\prod _{\goth p\in Z}\Hom _R(\LCMo ^i_{(y_1,\dots ,y_i)R}(R),\InjH
_R(R/\goth p))^{\mu_\goth p}$ is contained in some associated prime ideal of some $\Hom _R(\LCMo ^i_{(y_1,\dots ,y_i)R}(R),\InjH
_R(R/\goth p))$. But the set
$$\Ass _R(\Hom _R(\LCMo ^i_{(y_1,\dots ,y_i)R}(R),\InjH
_R(R/\goth p)))$$
has the ``openness property'' (O) because of the
conjecture, and so
$$\eqalign {\Ass _R(\Hom _R(\LCMo ^i_{(y_1,\dots ,y_i)R}(R),E))&=\Ass _R(\bigoplus _{\goth p\in Z}\Hom _R(\LCMo ^i_{(y_1,\dots ,y_i)R}(R),\InjH
_R(R/\goth p))^{(\mu_\goth p)})\cr &=\bigcup _{\goth p\in Z}\Ass _R(\Hom
_R(\LCMo ^i_{(y_1,\dots ,y_i)R}(R),\InjH _R(R/\goth p)))\cr &\subseteq \Ass
_R(D_R(\LCMo ^i_{(y_1,\dots ,y_i)R}(R)))\cr
}$$
In particular
$$\Ass _R(D_{R_0}(\LCMo ^i_{(y_1,\dots ,y_i)R}(R)))=\Ass _R(D_R(\LCMo
^i_{(y_1,\dots ,y_i)R}(R)))$$
\bigskip
\mittel
{\bf 2.4 An easy application}
\normal
\smallskip
The following result is well-known; we present a new proof for it:
\bigskip
{\bf 2.4.1 Corollary}
\par
Let $(R,\goth m)$ be a noetherian local ring and $M$ a
finitely generated $R$-module. Then $\LCMo ^{\dim
_R(M)}_\goth m(M)\neq 0$.
\par
Proof:
\par
We use the following fact: For every $n\in {\bf N} $, every ideal $I\subseteq
R$ and every finitely generated $R$-module $N$ the following statemtents are equivalent:
\par
(i) $H^i_I(N)=0$ for all $i\geq n$.
\par
(ii) $H^i_I(R/\goth p)=0$ for all $i\geq n$ and all $ \goth p\in
\Supp _R(N)$
\par
This fact implies (setting $d:=\dim _R(M)=\dim (R/\ann _R(M))$) $\LCMo ^d_\goth
m(M)\neq 0\iff \LCMo ^d_\goth m(R/\ann _R(M))$. Thus we may assume
$M=R$ and $R$ is a domain. Again we set $d:=\dim (R)$ and choose a system of parameters $x_1,\dots ,x_d\in R$
for $R$; theorems 2.2.1 and 2.2.3 (in the case of mixed characteristic we may
choose $x_1=p$ where $p=\char (R/\goth m)$) implie $\{ 0\} \in \Ass _R(D(\LCMo ^d_\goth m(R)))$;
in particular $\LCMo ^d_\goth m(R)\neq 0$.

\def\litem{\par\noindent \hangindent=\parindent\ltextindent}
\def\ltextindent#1{\hbox to \hangindent{#1\hss}\ignorespaces}
\bigskip
\mittel
{\bf References}
\normal
\smallskip
\parindent=0.8cm
\litem{1.} Bass, H. On the ubiquity of Gorenstein rings, {\it Math. Z.} {\bf
82}, (1963) 8-28.
\medskip
\litem{2.} Brodmann, M. and Hellus, M. Cohomological patterns of coherent
sheaves over projective schemes, {\it Journal of Pura and Applied Algebra}
{\bf 172}, (2002) 165-182.
\medskip
\litem{3.} Brodmann, M. P. and Sharp, R. J. Local Cohomology, {\it Cambridge
studies in advanced mathematics} {\bf 60}, (1998).
\medskip
\litem{4.} Bruns, W. and Herzog, J. Cohen-Macaulay Rings, {\it Cambridge
University Press}, (1993).
\medskip
\litem{5.} Eisenbud, D. Commutative Algebra with A View Toward Algebraic
Geometry, {\it SPringer Verlag}, (1995).
\medskip
\litem{6.} Grothendieck, A. Local Cohomology, {\it Lecture Notes in
Mathematics, Springer Verlag}, (1967).
\medskip
\litem{7.} Hellus, M. On the set of associated primes of a local cohomology
module, {\it J. Algebra} {\bf 237}, (2001) 406-419.
\medskip
\litem{8.} Huneke, C. Problems on Local Cohomology, {\it Res.
Notes Math. } {\bf 2}, (1992).
\medskip
\litem{9.} Matlis, E. Injective modules over Noetherian rings, {\it Pacific J. Math.} {\bf 8}, (1958) 511-528.
\medskip
\litem{10.} Matsumura, H. Commutative ring theory, {\it Cambridge
University Press}, (1986).
\medskip
\litem{11.} Scheja, G. and Storch, U. Regular Sequences and Resultants, {\it AK
Peters}, (2001).

\end